\documentclass[a4paper,12pt]{article}

\usepackage{amsmath,amssymb,latexsym,amsfonts,amsthm}
\usepackage{mathrsfs}

\theoremstyle{plain}
\newtheorem{Thm}{Theorem}[section]
\newtheorem{Lem}[Thm]{Lemma}
\newtheorem{Cor}[Thm]{Corollary}
\theoremstyle{definition}
\newtheorem{Rem}[Thm]{Remark}

\newcommand{\bN}{\ensuremath{\mathbb{N}}}
\newcommand{\bR}{\ensuremath{\mathbb{R}}}
\newcommand{\bT}{\ensuremath{\mathbb{T}}}
\newcommand{\bZ}{\ensuremath{\mathbb{Z}}}
\newcommand{\cB}{\ensuremath{\mathcal{B}}}
\newcommand{\cF}{\ensuremath{\mathcal{F}}}
\newcommand{\cG}{\ensuremath{\mathcal{G}}}
\newcommand{\fG}{\ensuremath{\mathfrak{G}}}
\newcommand{\fH}{\ensuremath{\mathfrak{H}}}
\newcommand{\sP}{\ensuremath{\mathscr{P}}}
\newcommand{\vh}{\ensuremath{\mbox{{\boldmath $h$}}}}
\newcommand{\vs}{\ensuremath{\mbox{{\boldmath $s$}}}}

\newcommand{\e}{{\rm e}}
\renewcommand{\d}{{\rm d}}

\newcommand{\dist}{\stackrel{\d}{=}}
\newcommand{\cdist}{\stackrel{\d}{\longrightarrow}}
\newcommand{\asto}{\stackrel{{\rm a.s.}}{\longrightarrow}}

\newcommand{\absol}[1]{\left| #1 \right|} 
\newcommand{\rbra}[1]{\left( #1 \right)} 
\newcommand{\cbra}[1]{\left\{ #1 \right\}} 
\newcommand{\sbra}[1]{\left[ #1 \right]} 

\renewcommand{\tilde}{\widetilde}

\numberwithin{equation}{section}

\makeatletter
\renewcommand\section{\@startsection {section}{1}{\z@}%
                                   {-3.5ex \@plus -1ex \@minus -.2ex}%
                                   {2.3ex \@plus.2ex}%
                                   {\normalfont\large\bf}}
\makeatother

\makeatletter
\@addtoreset{footnote}{page}
\makeatother

\begin{document}
\begin{center}
{\Large \bf Extremal solutions for stochastic equations 
indexed by negative integers 
and taking values in compact groups}
\end{center}
\begin{center}
Takao \textsc{Hirayama}\footnote{
Department of Mathematical Sciences,
Ritsumeikan University, Shiga, JAPAN.} 
and Kouji \textsc{Yano}\footnote{
Department of Mathematics, Graduate School of Science,
Kobe University, Kobe, JAPAN.}\footnote{
The research of this author was supported by KAKENHI (20740060)}
\end{center}
\begin{center}
{\small \today}               
\end{center}
\begin{center}
{\small Dedicated to Professor Marc Yor 
on the occasion of his 60th birthday}
\end{center}


\begin{abstract}
Stochastic equations indexed by negative integers 
and taking values in compact groups 
are studied. 
Extremal solutions of the equations are characterized 
in terms of infinite products of independent random variables. 
This result is applied to characterize several properties of the set of all solutions 
in terms of the law of the driving noise. 
\end{abstract}

\noindent
{\footnotesize Keywords and phrases: Stochastic equation, Markov process, extremal
point, infinite convolution, independent complement.} 
\\
{\footnotesize AMS 2000 subject classifications: 
Primary
60J05; 
secondary
60J50; 
60B15. 
} 

\section{Introduction}

Let $ G $ be a compact topological group. 
We consider the following stochastic equation on the state space $ G $ 
indexed by $ -\bN $: 
\begin{align}
\eta_k = \xi_k \eta_{k-1} 
, \quad k \in -\bN 
\label{eq: intro SE}
\end{align}
where $ (\eta_k) = (\eta_k: k \in -\bN) $ 
is an unknown process 
and $ (\xi_k) = (\xi_k: k \in -\bN) $ is a driving noise, 
i.e., the $ \xi_k $'s are independent (but in general not identically distributed), 
both taking values in $ G $. 
Iterating equation \eqref{eq: intro SE}, we have 
\begin{align}
\eta_k = \xi_k \xi_{k-1} \cdots \xi_{l+1} \eta_{l} 
, \quad k,l \in -\bN , \ k>l . 
\label{eq: etak and etak-n}
\end{align}
If we regard $ \eta_{l} $ as an initial state, 
then the states afterwards $ \{ \eta_{l+1},\eta_{l+2},\ldots,\eta_0 \} $ may be obtained from 
the noise $ \{ \xi_{l+1},\xi_{l+2},\ldots,\xi_0 \} $ 
together with the initial state $ \eta_l $. 
But the difficulty in the study of equation \eqref{eq: intro SE} 
comes from the fact that there is a priori 
no ``initial state at time $ -\infty $". 

We are interested in conditions on the noise law 
for the set of all possible solutions of equation \eqref{eq: intro SE} 
to satisfy certain properties. 
In particular, we raise the following three questions 
(all of which will be stated precisely in the next section): 
\\ \quad {\bf (Q1)} 
Does {\em uniqueness in law} hold? 
\\ \quad {\bf (Q2)} 
Does there exist a {\em strong solution}, 
i.e., a solution where each $ \eta_k $ is measurable 
with respect to the noise up to time $ k $? 
\\ \quad {\bf (Q3)} 
If a solution is non-strong, 
the noise process up to time $ k $ is inadequate 
to completely know the value of $ \eta_k $. 
Can we find some independent $ G $-valued random variable 
which complements the lack of information? 

The purpose of the present paper is 
to give clear answers to {\bf (Q1)}-{\bf (Q3)}. 
Our results generalize those of Yor \cite{Y} 
and complete those of Akahori et al.~\cite{AUY}. 
We point out that a key role is played by 
{\em extremal solutions}, 
which are precisely the solutions whose 
remote past is trivial (see Section \ref{sec: main}). 
For this purpose, 
we shall utilize 
the general theorems (see Theorems \ref{C2} and \ref{C3}) 
about infinite products of independent random variables, 
which are due to Kloss \cite{K}, Tortrat \cite{T} and Csisz{\'a}r \cite{Csi}. 
We will see that, thanks to the choice of $ -\bN $, instead of $ \bN $, 
as the index set, these theorems are deepened by 
our main theorem (Theorem \ref{thm: main1}) 
in terms of Markov processes. 

The present paper is organized as follows. 
In Section \ref{sec: main}, we give some notations, explain our terminology, 
and then state our main theorems. 
In Section \ref{sec: sigma fields}, 
we give two important lemmas concerning $ \sigma $-fields. 
In Section \ref{sec: backgr}, we recall some of the results 
of Yor \cite{Y} and Akahori et al.~\cite{AUY}. 
Section \ref{sec: infin prod} 
is devoted to the proofs of main theorems.

\section{Main results} \label{sec: main}

\subsection{Notations and terminology} \label{sec: term}

Let $ G $ be a compact topological group 
which we assume to be Hausdorff and with a countable basis. 
Then $ G $ is necessarily metrizable (see, e.g., \cite[Prop.7.1.12]{Cohn}); 
in particular, $ G $ is a Polish space. 
To avoid trivial complications, 
we suppose that $ G $ contains more than one element. 

Let us give precise definitions 
as to the terminology appearing in {\bf (Q1)}-{\bf (Q3)}, 
which is related to some filtration problems. 
Denote $ \bN = \{ 0,1,2,\ldots \} $. 
For two processes $ (\eta_k:k \in -\bN) $ and $ (\xi_k:k \in -\bN) $ 
defined on a common probability space, 
we consider the three filtrations: 
\begin{align}
\cF^{\eta}_k = \sigma(\eta_m:m \le k) 
, \quad 
\cF^{\xi}_k = \sigma(\xi_m:m \le k) 
, \quad \text{and} \quad 
\cF^{\eta,\xi}_k = \sigma(\eta_m,\xi_m:m \le k) . 
\label{}
\end{align}

Let $ \mu = (\mu_k: k \in -\bN) $ 
be a family of probability laws $ \mu_k $ on $ G $. 
By a {\em solution of equation \eqref{eq: intro SE} 
(with the noise law $ \mu $)}, 
we mean a pair of processes $ \{ (\eta_k),(\xi_k) \} $ 
defined on a probability space $ (\Omega,\cF,P) $ 
such that 
\begin{align}
\text{for any $ k \in -\bN $}, \quad 
\begin{cases}
\eta_k = \xi_k \eta_{k-1} \quad \text{a.s.}, \\
\text{$ \xi_k $ is independent of $ \cF^{\eta}_{k-1} $}, \\
\text{$ \xi_k $ has law $ \mu_k $}. 
\end{cases}
\label{eq: cond of solution}
\end{align}
This is equivalent to stating that 
$ (\eta_k) $ is a (possibly time-inhomogeneous) Markov process indexed by $ -\bN $ 
such that 
\begin{align}
E[f(\eta_k) | \cF^{\eta}_{k-1}] 
= \int_G f(g \eta_{k-1}) \mu_k(\d g) 
\quad \text{a.s.} 
, \quad k \in -\bN 
\label{eq: Markov}
\end{align}
for all non-negative Borel function $ f $ on $ G $. 
We note that, since 
\begin{align}
\xi_k = \eta_k (\eta_{k-1})^{-1} 
\quad \text{a.s.}, 
\label{eq: noise subset eta}
\end{align}
there is the equality $ \cF^{\eta,\xi}_k = \cF^{\eta}_k $; 
in particular, $ \cF^{\xi}_k \subset \cF^{\eta}_k $. 

Let us fix $ \mu = (\mu_k: k \in -\bN) $ throughout this paper. 
Following Yor \cite{Y} and Akahori et al. \cite{AUY}, 
we introduce the following definitions: 
\begin{itemize}
\item 
Let $ \sP_{\mu} $ denote the set of the laws of $ (\eta_k) $ on $ G^{-\bN} $ 
for all possible solutions $ \{ (\eta_k),(\xi_k) \} $. 
The set $ \sP_{\mu} $ is a compact convex subset of $ \sP(G^{-\bN}) $, 
the set of all probability laws on $ G^{-\bN} $, 
equipped with the topology of weak convergence 
and with the usual convex structure. 
\item 
We say that a solution $ \{ (\eta^0_k),(\xi_k) \} $ is {\em extremal} if 
the law of $ (\eta^0_k) $ is an extremal point of the compact convex set $ \sP_{\mu} $. 
We denote by $ \sP_{\mu}^{\rm extremal} $ 
the set of all extremal points of $ \sP_{\mu} $. 
\item 
We say that {\em uniqueness in law} holds 
if any two solutions have the same laws. 
\item 
We say that a solution $ \{ (\eta_k),(\xi_k) \} $ is {\em strong} if 
each $ \eta_k $ is a.s. measurable with respect to the past noise, i.e., 
$ \cF^{\eta}_k \subset \cF^{\xi}_k $ a.s. for all $ k \in - \bN $; 
so that in this case $ \cF^{\eta}_k = \cF^{\xi}_k $ a.s. 
We denote by $ \sP_{\mu}^{\rm strong} $ 
the set of the laws of $ (\eta_k) $ 
for all strong solutions $ \{ (\eta_k),(\xi_k) \} $. 
\end{itemize}

The terms ``uniqueness in law" and ``strong" 
originate from the theory of stochastic differential equations; 
see, e.g., \cite{IW}.

\subsection{Basic facts}

Let us recall some basic facts concerning solutions of equation \eqref{eq: intro SE}. 

First, we state without proof the following five facts numbered from {\bf 1)} to {\bf 5)}, 
which are due to Yor \cite{Y} and Akahori et al. \cite{AUY}. 
We will give their proofs in Section \ref{sec: backgr} for completeness of this paper. 

{\bf 1).} 
For a solution $ \{ (\eta_k),(\xi_k) \} $, 
the joint law of $ ((\eta_k),(\xi_k)) $ on $ G^{-\bN} \times G^{-\bN} $ 
is determined from the sequence $ (\lambda_k) $ 
of the marginal laws of $ (\eta_k) $ on $ G^{-\bN} $. 
We can and do in what follows identify two solutions having common joint laws, 
so that a solution $ \{ (\eta_k),(\xi_k) \} $ will be identified with 
the sequence $ (\lambda_k) $ 
as well as with the law of $ (\eta_k) $ on $ G^{-\bN} $, 
which is a point of $ \sP_{\mu} $. 

{\bf 2).} 
For a solution $ \{ (\eta_k),(\xi_k) \} $, 
the sequence $ (\lambda_k) $ of the marginal laws of $ (\eta_k) $ 
satisfies the convolution equation 
\begin{align}
\lambda_k = \mu_k * \lambda_{k-1} 
, \quad k \in -\bN . 
\label{eq: conv eq}
\end{align}
Conversely, for a sequence $ (\lambda_k) \subset \sP(G) $ satisfying 
the convolution equation \eqref{eq: conv eq}, 
there exists a solution $ \{ (\eta_k),(\xi_k) \} $ 
whose joint law is unique 
such that $ (\lambda_k) $ is the marginal laws of $ (\eta_k) $, 
i.e., $ \lambda_k = P(\eta_k \in \cdot) $ for each $ k \in -\bN $. 

{\bf 3).} 
Whatever the noise law $ \mu=(\mu_k) $ is, 
there always exists a solution $ \{ (\eta^*_k),(\xi_k) \} $ 
such that each $ \eta^*_k $ is uniform on $ G $, 
i.e., the law of $ \eta^*_k $ on $ G $ is the {\em normalized Haar measure} of $ G $. 
This is the case because $ G $ is compact. 
We call $ \{ (\eta^*_k),(\xi_k) \} $ {\em the uniform solution} 
and we denote the law of $ (\eta^*_k) $ by $ P^*_{\mu} $. 
From this, we obtain the following: 
\begin{itemize}
\item 
Uniqueness in law holds if and only if $ \sP_{\mu} = \{ P^*_{\mu} \} $. 
\item 
Since $ \sP_{\mu} $ is non-empty, so is $ \sP_{\mu}^{\rm extremal} $. 
\end{itemize}
Moreover, the uniform solution is non-strong; 
in fact, each $ \eta^*_k $ is independent of $ \cF^{\xi}_0 $. 

{\bf 4).} 
A solution $ \{ (\eta^0_k),(\xi_k) \} $ is extremal if and only if 
the remote past $ \cF^{\eta^0}_{-\infty } := \cap_k \cF^{\eta^0}_k $ is trivial. 
By Kolmogorov's 0-1 law, we see that 
a strong solution is always extremal; in other words, 
\begin{align}
\sP_{\mu}^{\rm strong} \subset \sP_{\mu}^{\rm extremal} \subset \sP_{\mu} . 
\label{}
\end{align}

{\bf 5).} 
Let an extremal solution $ \{ (\eta^0_k),(\xi_k) \} $ be fixed. 
Any other extremal solution is then identical in law to 
$ \{ (\eta^0_k g),(\xi_k) \} $ for some $ g \in G $. 
This shows that 
any solution is identical in law to 
$ \{ (\eta^0_k V),(\xi_k) \} $ 
for some $ G $-valued random variable $ V $ 
independent of $ \{ (\eta^0_k),(\xi_k) \} $. 

Second, we mention the following trichotomy, 
which may be deduced immediately from the above facts {\bf 1)}-{\bf 5)} 
(see also \cite{YY}): 
\\ {\bf Case (A): 
Uniqueness in law holds}, i.e., $ \sP_{\mu} = \{ P^*_{\mu} \} $. 
In this case, the uniform solution is the only solution, 
so it is extremal, 
but it is non-strong. 
\\ {\bf Case (B): 
There exists a strong solution}, i.e., $ \sP_{\mu}^{\rm strong} \neq \emptyset $. 
In this case, uniqueness in law fails. 
Moreover, it holds that $ \sP_{\mu}^{\rm strong} = \sP_{\mu}^{\rm extremal} $, 
i.e., all extremal solutions are strong, 
and all non-extremal ones are not. 
\\ {\bf Case (C): 
Uniqueness in law fails and there is no strong solution.} 
In this case, it holds that $ \sP_{\mu}^{\rm strong} = \emptyset $ 
and $ \sP_{\mu}^{\rm extremal} \subsetneq \sP_{\mu} $.

Third, we discuss some problem of filtrations. 
The discussion in the following seems elementary 
but needs more carefulness than one may expect; 
see, e.g., \cite[\S 2.5]{CY} and references therein. 
In fact, it has been a source of errors; 
see, e.g., \cite[a) of \S 5]{YY}. 

For decreasing $ \sigma $-fields $ \cF_0,\cF_{-1},\ldots $ 
and a $ \sigma $-field $ \cG $, 
it is obvious that 
\begin{align}
\bigcap_{l \in -\bN} ( \cF_{l} \vee \cG ) 
\supset 
\rbra{ \bigcap_{l \in -\bN} \cF_{l} } \vee \cG . 
\label{eq: simga fields inclusion2}
\end{align}
But this inclusion is sometimes strict; 
one cannot always change the order of the two operations $ \cap_{l \in -\bN} $ and $ \vee $. 
We will give in Lemma \ref{lem: interchange} some sufficient condition 
so that the equality holds in \eqref{eq: simga fields inclusion2}. 
For some recent discussions of this well-studied question, 
see Crimaldi et al.~\cite{CLP} and Berti et al.~\cite{BPR}. 

Let $ \{ (\eta_k),(\xi_k) \} $ be a solution. 
By equation \eqref{eq: intro SE}, we have 
\begin{align}
\bigcap_{l \in -\bN} ( \cF^{\eta}_{l} \vee \cF^{\xi}_0 ) = \cF^{\eta}_0 
\label{eq: sigma field A}
\end{align}
whereas 
\begin{align}
\rbra{ \bigcap_{l \in -\bN} \cF^{\eta}_{l} } \vee \cF^{\xi}_0 
= \cF^{\eta}_{-\infty } \vee \cF^{\xi}_0 . 
\label{eq: sigma field B}
\end{align}
As we have noted above, the $ \sigma $-field \eqref{eq: sigma field B} 
may be strictly contained 
in \eqref{eq: sigma field A}; 
in other words, 
the present $ \cF^{\eta}_0 $ may possess some extra information 
which cannot be explained 
by the noise $ \cF^{\xi}_0 $ together with $ \cF^{\eta}_{-\infty } $, 
the ``initial state at $ -\infty $". 
So we want to find a sub $ \sigma $-field $ \cG $ 
such that 
\begin{align}
\cF^{\eta}_0 =& \cG \vee \cF^{\eta}_{-\infty } \vee \cF^{\xi}_0 
\ \text{a.s.} 
\label{} \\
\text{and} \ & \text{the three $ \sigma $-fields 
$ \cG $, $ \cF^{\eta}_{-\infty } $ and $ \cF^{\xi}_0 $ are independent} . 
\label{}
\end{align}
We call $ \cG $ an {\em independent complement} 
of $ \cF^{\eta}_{-\infty } \vee \cF^{\xi}_0 $ in $ \cF^{\eta}_0 $. 
See Chaumont--Yor \cite[\S 2]{CY} and references therein.

\subsection{Extremal solution} 

Let us present our main theorems. 
The proofs of all theorems and corollaries presented in this subsection 
will be given in Section \ref{sec: infin prod}. 

For a compact subgroup $ H $ of $ G $, 
we denote by $ \omega_H $ the normalized Haar measure on $ H $. 
We denote by $ G/H $ the quotient set, i.e., the set 
of all left cosets $ gH = \{ gh:h \in H \} $ for $ g \in G $. 
The set $ G/H $ is equipped with the smallest topology 
in which the canonical projection $ G \ni g \mapsto gH \in G/H $ is continuous. 
Then we see that $ G/H $ is compact and metrizable. 

The following theorem, 
which is essentially due to Csisz\'ar \cite{Csi}, 
concerns infinite convolution products of probability laws on $ G $. 

\begin{Thm} \label{thm: main0}
There exist a sequence $ (\lambda_k) $ 
of probability laws on $ G $, 
a sequence $ (\alpha _l) $ of deterministic elements of $ G $, 
and a compact subgroup $ H $ of $ G $, such that 
the following statements hold: 
\\ \quad {\bf (I1)} 
$ \mu_k * \mu_{k-1} * \cdots * \mu_{l} * \delta_{\alpha _l} 
\to \lambda_k $ as $ l \to -\infty $ for each $ k \in -\bN $; $ \Big. $ 
\\ \quad {\bf (I2)} 
$ \delta_{\alpha _l^{-1}} * \lambda_{l-1} \to \omega_H $ 
as $ l \to -\infty $; $ \Big. $ 
\\ \quad {\bf (I3)} 
$ \lambda_k * \delta_h = \lambda_k $ 
for each $ h \in H $ and each $ k \in -\bN $; $ \Big. $ 
\\ \quad {\bf (I4)} 
If $ \xi_k $'s are independent random variables 
such that each $ \xi_k $ has law $ \mu_k $, $ \Big. $ 
\\ \quad 
then, for any $ k \in -\bN $, 
$ \xi_k \xi_{k-1} \cdots \xi_l \alpha _l H $ 
converges a.s. in $ G/H $ as $ l \to - \infty $. $ \Big. $ 
\\
If, moreover, 
$ (\tilde{\lambda}_k) $, $ (\tilde{\alpha }_k) $ and $ \tilde{H} $ 
also satisfy {\bf (I1)}-{\bf (I4)}, 
then it holds that 
\begin{align}
\tilde{\lambda}_k = \lambda_k * \delta_g 
, \quad 
\tilde{H} = g^{-1} H g 
\label{}
\end{align}
for all accumulation point $ g $ 
of $ \{ \alpha _{l}^{-1} \tilde{\alpha }_{l} : l \in -\bN \} $. 
\end{Thm}

We remark that the sequence $ (\lambda_k) $ above 
satisfies the convolution equation \eqref{eq: conv eq}. 
This suggests that it is natural 
to choose $ -\bN $, instead of $ \bN $, 
as the index set. 
We may say that 
the following theorem, which characterizes extremal solutions, 
deepens Theorem \ref{thm: main0} in terms of Markov processes. 

\begin{Thm} \label{thm: main1}
For any extremal solution $ \{ (\eta^0_k),(\xi_k) \} $, 
there exist a sequence $ (\alpha _l) $ of deterministic elements of $ G $ 
and a compact subgroup $ H $ of $ G $ 
such that the following four conditions hold: 
\\ \quad {\bf (E1)} 
$ \xi_k \xi_{k-1} \cdots \xi_{l} \alpha _{l} \cdist \eta^0_k $ 
as $ l \to -\infty $ for each $ k \in -\bN $; $ \Big. $ 
\\ \quad {\bf (E2)} 
$ \alpha _l^{-1} \eta^0_{l-1} \cdist U_H $ 
as $ l \to -\infty $ 
where $ U_H $ is uniform on $ H $; $ \Big. $ 
\\ \quad {\bf (E3)} 
$ (\eta^0_k h) \dist (\eta^0_k) $ 
for each $ h \in H $ and each $ k \in -\bN $; $ \Big. $ 
\\ \quad {\bf (E4)} 
$ \xi_k \xi_{k-1} \cdots \xi_{l} \alpha _{l} H \asto \eta^0_k H $ 
as $ l \to -\infty $ 
for each $ k \in -\bN $. $ \Big. $ 
\\
If, moreover, $ \{ (\tilde{\eta}^0_k),(\tilde{\xi}_k) \} $ 
is another extremal solution which satisfies {\bf (E1)}-{\bf (E4)} 
with $ (\tilde{\alpha }_l) $ and $ \tilde{H} $, 
then it holds that 
\begin{align}
(\tilde{\eta}^0_k) \dist (\eta^0_k g) 
, \quad 
\tilde{H} = g^{-1} H g 
\label{}
\end{align}
for all accumulation point $ g $ 
of $ \{ \alpha _{l}^{-1} \tilde{\alpha }_{l} : l \in -\bN \} $. 
\end{Thm}

\begin{Rem}
Since $ \eta^0_{k-1} $ 
is independent of $ \sigma(\xi_k, \xi_{k-1}, \ldots \xi_{l}) $, 
we can combine the two conditions {\bf (E1)} and {\bf (E2)} together as follows: 
\\ \quad $ \bullet $ 
$ (\xi_k \xi_{k-1} \cdots \xi_{l} \alpha _{l},\alpha _l^{-1} \eta^0_{l-1}) 
\cdist (\eta^0_k,U_H) $ 
as $ l \to -\infty $ for each $ k \in -\bN $ $ \Big. $ 
\\ \quad 
where $ U_H $ is independent of $ \eta^0_k $ and is uniform on $ H $. 
\\
In particular, this shows $ \eta^0_k U_H \dist \eta^0_k $, 
which immediately implies {\bf (E3)}. 
\end{Rem}

The following theorem plays an essential role in the proof of Theorem \ref{thm: main1}. 

\begin{Thm} \label{thm: main1.2}
The following statements hold: 
\\ {\rm (i)} 
If $ \{ (\eta^0_k),(\xi_k) \} $ is a solution 
and satisfies {\bf (E1)} with some $ (\alpha _l) $, 
then it is extremal. 
\\ {\rm (ii)} 
If $ \{ (\eta^0_k),(\xi_k) \} $ is a solution 
which satisfies {\bf (E1)}, {\bf (E3)} and {\bf (E4)} 
for some $ (\alpha _l) $ and some $ H $, 
then it also satisfies {\bf (E2)} with these $ (\alpha _l) $ and $ H $. 
\end{Thm}

For a given noise law $ \mu=(\mu_k) $, 
the compact subgroup $ H $ of Theorem \ref{thm: main0} (or Theorem \ref{thm: main1}) 
is unique up to conjugacy, 
so we sometimes denote it by $ H_{\mu} $. 
The subgroup $ H_{\mu} $ may be characterized as follows. 

\begin{Cor} \label{thm: main1.5}
Let $ \{ (\eta^0_k),(\xi_k) \} $, $ (\alpha _l) $ and $ H_{\mu} $ 
be as in Theorem \ref{thm: main1}. 
Then the following statements hold: 
\\ \quad {\rm (i)} 
$ H_{\mu} = \cbra{ h \in G : (\eta^0_k h) \dist (\eta^0_k) } $; 
\\ \quad {\rm (ii)} 
$ H_{\mu} $ is the smallest compact subgroup $ H $ such that 
$ \xi_k \xi_{k-1} \cdots \xi_{l} \alpha _{l} H \asto \eta^0_k H $ 
as $ l \to -\infty $ for all $ k \in -\bN $. 
\end{Cor}

The following corollary answers {\bf (Q1)} and {\bf (Q2)}. 

\begin{Cor} \label{thm: main2}
Let $ H_{\mu} $ be as in Theorem \ref{thm: main0}. 
Then the following statements hold: 
\\ {\rm (A)} 
The following statements are equivalent: 
\\ \quad {\rm (A1)} 
Uniqueness in law holds; 
\\ \quad {\rm (A2)} 
$ H_{\mu} = G $; 
\\ \quad {\rm (A3)} 
For any $ k \in -\bN $, 
$ \xi_k \xi_{k-1} \cdots \xi_{l} \cdist U_G $ as $ l \to -\infty $ 
where $ U_G $ is uniform on $ G $. 
\\ {\rm (B)} 
The following statements are equivalent: 
\\ \quad {\rm (B1)} 
There exists a strong solution; 
\\ \quad {\rm (B2)} 
$ H_{\mu} = \{ {\rm unit} \} $; 
\\ \quad {\rm (B3)} 
There exists a sequence $ (\alpha _l) $ of deterministic elements of $ G $ 
such that, for each $ k \in -\bN $, $ \xi_k \xi_{k-1} \cdots \xi_{l} \alpha _{l} $ 
converges a.s. as $ l \to -\infty $. 
\\
In this case, if we write $ \eta^0_k $ for the limit in {\rm (B3)}, 
then the pair $ \{ (\eta^0_k),(\xi_k) \} $ is 
a strong (and consequently extremal) solution. 
\end{Cor}

\subsection{Complementation formulae} \label{sec: factor and compl}

For any compact subgroup $ H $ of $ G $, 
there always exists a measurable section $ \vs(\cdot):G/H \to G $ 
(see \cite[Exercise 8.4]{Cohn}). 
We define the measurable mapping $ \vh(\cdot):G \to H $ as 
\begin{align}
\vh(g) = (\vs(gH))^{-1} g 
, \quad g \in G . 
\label{eq: measurable decomp}
\end{align}
Then the mapping 
\begin{align}
G \ni g = \vs(gH) \vh(g) \mapsto (gH,\vh(g)) \in (G/H) \times H 
\label{eq: group decomp}
\end{align}
is a bi-measurable bijection, 
where the direct product $ (G/H) \times H $ 
is equipped with product topology. 

In this subsection, 
we assume that $ \mu=(\mu_k) $ denotes a given noise law 
and that $ (\lambda_k) $, $ (\alpha _l) $ and $ H $ are 
as in Theorem \ref{thm: main0}. 
Let $ \vs(\cdot):G/H \to G $ be a measurable section 
associated with this $ H $ 
and $ \vh(\cdot) $ be as defined by \eqref{eq: measurable decomp}. 

The following theorem provides us with a procedure 
of constructing an extremal solution 
from the noise together with an additional randomness.

\begin{Thm} \label{thm: main1.4}
Let $ \xi_k $'s be independent random variables such that each $ \xi_k $ has law $ \mu_k $ 
and let $ U_0 $ be a $ G $-valued random variable independent of $ (\xi_k) $. 
For each $ k \in -\bN $, define 
\begin{align}
\phi_k =& \vs \rbra{ \lim_{l \to -\infty } \xi_k \xi_{k-1} \cdots \xi_l \alpha _l H } , 
\label{} \\
U_k =& \phi_k^{-1} (\xi_0 \xi_{-1} \cdots \xi_{k+1})^{-1} \phi_0 U_0 , 
\label{}
\end{align}
and then define 
\begin{align}
\eta^0_k =& \phi_k U_k . 
\label{}
\end{align}
Then $ \{ (\eta^0_k),(\xi_k) \} $ is an extremal solution 
such that each $ \eta^0_k $ has law $ \lambda_k $. 
Moreover, for any $ k \in -\bN $, it holds that 
$ U_k $ is independent of $ \cF^{\xi}_0 $ and is uniform on $ H $, 
and that 
\begin{align}
\cF^{\eta^0}_k = \sigma(U_k) \vee \cF^{\xi}_k 
\label{}
\end{align}
where $ \sigma(U_k) $ and $ \cF^{\xi}_k $ are independent. 
\end{Thm}

Theorem \ref{thm: main1.4} will be proved in Subsection \ref{sec: construction}.

By Theorem \ref{thm: main1.4} and by point {\bf 4)} of Subsection \ref{sec: term}, 
any solution $ \{ (\eta_k),(\xi_k) \} $ may be represented as 
$ \eta_k = \phi_k U_k V $ 
for some random variable $ V $ independent of $ \cF^{\xi}_0 \vee \sigma(U_0) $, 
and consequently, it holds that, for any $ k \in -\bN $, 
\begin{align}
\cF^{\eta}_k \subset \sigma(U_k) \vee \sigma(V) \vee \cF^{\xi}_k 
\quad \text{a.s.} 
\label{}
\end{align}
For the converse inclusion, we need to take $ V $ nicely 
and to represent $ V $ and $ U_k $ in terms of $ (\eta_k) $. 
The following theorem solves this problem 
and answers {\bf (Q3)} completely. 

\begin{Thm} \label{thm: main3}
Let $ \{ (\eta_k),(\xi_k) \} $ be any solution. 
For any $ k \in -\bN $, 
define 
\begin{align}
\phi_k =& 
\vs \rbra{ \lim_{l \to -\infty } \xi_k \xi_{k-1} \cdots \xi_l \alpha _l H } , 
\label{eq: main4 01} \\
V =& \vs \rbra{ \lim_{l \to -\infty } \eta_l^{-1} \phi_l H }^{-1} , 
\label{eq: main4 02} \\
U_k =& \vh(\eta_k V^{-1}) . 
\label{eq: main4 03}
\end{align}
Then, for any $ k \in -\bN $, 
the random variable $ \eta_k $ is factorized as $ \eta_k = \phi_k U_k V $ 
and the following statements hold: 
\\ \quad {\rm (i)} 
$ \phi_k \in \cF^{\xi}_k $ a.s.; 
\\ \quad {\rm (ii)} 
$ U_k $ is independent of $ \sigma(V) \vee \cF^{\xi}_0 $ and is uniform on $ H $; 
\\ \quad {\rm (iii)} 
$ \cF^{\eta}_{-\infty } = \sigma(V) $. 
\\
Moreover, it holds that, for any $ k \in -\bN $, 
\begin{align}
\cF^{\eta}_k = \sigma(U_k) \vee \sigma(V) \vee \cF^{\xi}_k 
\quad \text{a.s.} 
\label{eq: main3 sigma fields}
\end{align}
where the three $ \sigma $-fields 
$ \sigma(U_k) $, $ \sigma(V) $ and $ \cF^{\xi}_k $ 
are independent. 
\end{Thm}

Theorem \ref{thm: main3} will be proved in Subsection \ref{sec: compl}. 

\begin{Cor} \label{thm: main4}
Let $ \{ (\eta_k),(\xi_k) \} $ be any solution. 
Then the identity 
\begin{align}
\bigcap_{l \in -\bN} ( \cF^{\eta}_{l} \vee \cF^{\xi}_0 ) 
= 
\rbra{ \bigcap_{l \in -\bN} \cF^{\eta}_{l} } \vee \cF^{\xi}_0 
\quad \text{a.s.} 
\label{eq: sigma field C}
\end{align}
holds 
if and only if there exists a strong solution. 
\end{Cor}

\begin{proof}
By \eqref{eq: sigma field A} and \eqref{eq: sigma field B}, 
the identity \eqref{eq: sigma field C} holds if and only if 
$ \cF^{\eta}_0 = \cF^{\eta}_{-\infty } \vee \cF^{\xi}_0 $. 
By Theorem \ref{thm: main3}, this is equivalent to 
triviality of $ \sigma(U_0) $, 
which leads to $ H=\{ {\rm unit} \} $. 
The proof is now completed by (B) of Corollary \ref{thm: main2}. 
\end{proof}

\subsection{The case of one-dimensional torus}

Let us consider the case of one-dimensional torus $ G = \bT \cong [0,1) $. 
In this case we prefer addition instead of multiplication, 
so that equation \eqref{eq: intro SE} may be rewritten as 
\begin{align}
\eta_k = \xi_k + \eta_{k-1} 
, \quad k \in -\bN . 
\label{}
\end{align}
For a given noise law $ \mu=(\mu_k) $, 
the compact subgroup $ H=H_{\mu} $ as in Theorem \ref{thm: main0} 
is uniquely determined. 
We have the following three distinct cases: 
\\ \quad {\rm (A)} 
$ H_{\mu}=[0,1) $. 
\\ \quad {\rm (B)} 
$ H_{\mu}=\{ 0 \} $. 
\\ \quad {\rm (C)} 
$ H_{\mu} $ may be expressed as 
\begin{align}
H_{\mu} = \cbra{ 0, \frac{1}{p_{\mu}}, \ldots, \frac{p_{\mu}-1}{p_{\mu}} } 
\label{eq: H mu}
\end{align}
for some integer $ p_{\mu} \ge 2 $. 

For $ x \in \bR $, we write $ [x] $ for the integer part of $ x $, 
i.e., the largest integer which does not exceed $ x $, 
and write $ \{ x \} $ for the fractional part of $ x $, 
i.e., $ \{ x \} = x-[x] $. 
In the case {\rm (C)}, 
we identify the quotient set $ G/H_{\mu} $ with $ [0,1/p_{\mu}) $ ($ \cong \bT $). 
In this case, we may choose as the measurable section $ \vs(\cdot) $ the mapping 
\begin{align}
\vs_{\mu}(x+H_{\mu}) = \{ p_{\mu} x \} / p_{\mu} 
, \quad x \in [0,1) , 
\label{}
\end{align}
hence we see that $ \vh(\cdot) = \vh_{\mu}(\cdot) $ is given as 
\begin{align}
\vh_{\mu}(x) = [p_{\mu} x] / p_{\mu} 
, \quad x \in [0,1) . 
\label{}
\end{align}

Now we obtain the following corollary (see Subsection \ref{sec: Yor}): 

\begin{Cor} \label{thm: main5}
Suppose that $ G = \bT \cong [0,1) $. 
Let $ (\alpha _l) $ and $ H_{\mu} $ 
as in Theorem \ref{thm: main1}. 
Let $ \{ (\eta_k),(\xi_k) \} $ be any solution. 
Then the following statements hold: 
\\ \quad {\rm (A)} 
Uniqueness in law holds if and only if $ H_{\mu}=[0,1) $. 
In this case, it holds that $ \cF^{\eta}_{-\infty } $ is trivial and, 
for any $ k \in -\bN $, that 
$ \eta_k $ is uniform on $ G $ and 
\begin{align}
\cF^{\eta}_k = \sigma(\eta_k) \vee \cF^{\xi}_k 
\quad \text{a.s.} 
\label{}
\end{align}
where $ \sigma(\eta_k) $ and $ \cF^{\xi}_k $ are independent; 
\\ \quad {\rm (B)} 
There exists a strong solution if and only if $ H_{\mu} = \{ 0 \} $. 
In this case, for any $ k \in -\bN $, the limits 
\begin{align}
\phi_k := \lim_{l \to -\infty } \rbra{ \sum_{j=l}^k \xi_j + \alpha _l } 
, \quad 
V := \lim_{l \to -\infty } (\eta_l - \phi_l) 
\quad \text{a.s.}, 
\label{}
\end{align}
exist and $ \eta_k = \phi_k + V $ 
where $ V $ is independent of $ \cF^{\xi}_0 $. 
Moreover, it holds that $ \cF^{\eta}_{-\infty } = \sigma(V) $ a.s. and, 
for any $ k \in -\bN $, that 
\begin{align}
\cF^{\eta}_k = \sigma(V) \vee \cF^{\xi}_k 
\quad \text{a.s.} 
\label{}
\end{align}
where $ \sigma(V) $ and $ \cF^{\xi}_k $ are independent; 
\\ \quad {\rm (C)} 
Suppose that $ H_{\mu} $ is of the form \eqref{eq: H mu} 
for $ p_{\mu} \ge 2 $. 
Then, for any $ k \in -\bN $, the limits 
\begin{align}
\phi_k :=& \lim_{l \to -\infty } 
\rbra{ \sum_{j=l}^k \xi_j + \alpha _l } 
\ \text{modulo} \ \frac{1}{p_{\mu}} 
\quad \text{a.s.}, 
\label{} \\
V :=& \lim_{l \to -\infty } \rbra{  \eta_l - \phi_l } 
\ \text{modulo} \ \frac{1}{p_{\mu}} 
\quad \text{a.s.} 
\label{}
\end{align}
exist where $ V $ is independent of $ \cF^{\xi}_0 $; 
\begin{align}
U_k := \frac{1}{p_{\mu}} \sbra{ p_{\mu} \rbra{ \eta_k - V } } 
\label{}
\end{align}
is uniform on $ H_{\mu} $ and is independent of $ \sigma(V) \vee \cF^{\xi}_0 $. 
Moreover, for any $ k \in -\bN $, it holds that 
$ \eta_k = \phi_k + U_k + V $ 
and that 
\begin{align}
\cF^{\eta}_k = \sigma(U_k) \vee \sigma(V) \vee \cF^{\xi}_k 
\quad \text{a.s.} 
\label{}
\end{align}
where 
the three $ \sigma $-fields 
$ \sigma(U_k) $, $ \sigma(V) $ and $ \cF^{\xi}_k $ 
are independent. 
\end{Cor}

\section{Some discussion on $ \sigma $-fields} \label{sec: sigma fields}

In this section 
we give two lemmas concerning $ \sigma $-fields, 
which will play important roles in our analysis. 
These lemmas seem elementary but should be dealt with carefully, 
because their statements are sources of errors. 
The first one is as follows. 

\begin{Lem} \label{lem: independence}
Let $ (\Omega,\cF,P) $ be a probability space. 
Let $ \cF_0 $ and $ \cG $ be two sub $ \sigma $-fields. 
Let $ X $ be an integrable random variable. 
Suppose that $ \sigma(X) \vee \cF_0 $ is independent of $ \cG $. 
Then it holds that 
\begin{align}
E[X|\cF_0 \vee \cG] = E[X|\cF_0] . 
\label{eq: independence conclusion}
\end{align}
\end{Lem}

\begin{proof}
Since $ \sigma(X) \vee \cF_0 $ is independent of $ \cG $, 
we have, for $ A \in \cF_0 $ and $ B \in \cG $, 
\begin{align}
E[X 1_A 1_B] = E[X 1_A] E[1_B] 
= E[E[X|\cF] 1_A] E[1_B] 
= E[E[X|\cF] 1_A 1_B] . 
\label{}
\end{align}
Thus, a monotone class argument yields 
\begin{align}
E[X 1_C] = E[E[X|\cF_0] 1_C] 
, \quad C \in \cF_0 \vee \cG . 
\label{}
\end{align}
Now the proof is complete. 
\end{proof}

\begin{Rem}
If we assume, instead of the independence between $ \sigma(X) \vee \cF_0 $ and $ \cG $, 
that $ \cF_0 $ and $ \cG $ are independent, 
then the conclusion \eqref{eq: independence conclusion} does not hold; 
see \cite[Exercise 2.2.1]{CY} for counterexamples. 
\end{Rem}

The second one is taken from \cite[Exercise 2.5.1]{CY}. 

\begin{Lem} \label{lem: interchange}
Let $ (\Omega,\cF,P) $ be a probability space. 
Let $ \{ \cF_{0}, \ \cF_{-1} , \ldots \} $ 
be a decreasing sequence of sub $ \sigma $-fields 
and $ \cG $ a sub $ \sigma $-field. 
Suppose that $ \cF_{0} $ is independent of $ \cG $. 
Then it holds that 
\begin{align}
\bigcap_{l \in -\bN } \rbra{ \cF_{l} \vee \cG } = 
\rbra{ \bigcap_{l \in -\bN } \cF_{l} } \vee \cG . 
\label{}
\end{align}
\end{Lem}

\begin{proof}
Let us write $ \cF_{-\infty } $ for $ \cap_{l \in -\bN } \cF_{l} $. 
It suffices to prove that 
$ \cap_{l \in -\bN } ( \cF_{l} \vee \cG ) \subset \cF_{-\infty } \vee \cG $, 
since the opposite inclusion is obvious. 
Let $ A \in \cF_{0} $ and $ B \in \cG $. 
Then, on one hand, we have 
\begin{align}
E[1_A 1_B|\cF_{l} \vee \cG] 
\stackrel{l \to -\infty }{\longrightarrow} 
E \sbra{ 1_A 1_B \bigg| \bigcap_{l \in -\bN } (\cF_{l} \vee \cG) } . 
\label{}
\end{align}
On the other hand, we have 
\begin{align}
E[1_A 1_B|\cF_{l} \vee \cG] 
=& E[1_A |\cF_{l} \vee \cG] 1_B 
= E[1_A |\cF_{l}] 1_B 
\quad \text{(from Lemma \ref{lem: independence})} 
\\
\stackrel{l \to -\infty }{\longrightarrow}& 
E[1_A |\cF_{-\infty }] 1_B 
= E[1_A |\cF_{-\infty } \vee \cG] 1_B 
= E[1_A 1_B |\cF_{-\infty } \vee \cG] . 
\label{}
\end{align}
Hence we see that the identity 
\begin{align}
E \sbra{ X \bigg| \bigcap_{l \in -\bN } (\cF_{l} \vee \cG) } = 
E [ X | \cF_{-\infty } \vee \cG] 
\label{eq: interchange (1)}
\end{align}
holds for $ X = 1_{A \cap B} $. 
A monotone class argument shows that the identity \eqref{eq: interchange (1)} 
holds for all $ X \in L^1(\cF_{0} \vee \cG) $. 
Now the proof is complete. 
\end{proof}

\section{Stochastic equations indexed by negative integers} \label{sec: backgr}

Our problem originates from 
Tsirelson's example of a stochastic differential equation with driving Brownian motion 
which has no strong solution (\cite{C}). 
He reduced the problem 
to equation \eqref{eq: intro SE} on the torus $ G = \bR / \bZ $ 
where the noise process consists of the projections of independent Gaussian variables. 
See \cite{YT} and \cite{YY} for brief surveys of this topic; 
see also \cite{Takahashi}. 
(Note that, in \cite{YT} and \cite{Takahashi}, 
the authors used the word ``remote past" 
for ``something at the time $ -\infty $", 
which is misleading 
because it is different from the usual terminology 
where ``remote past" means the $ \sigma $-field $ \cF^{\eta}_{-\infty } $.)

\subsection{Yor's stochastic equation} \label{sec: Yor}

Looking for some better understanding of the properties of Tsirelson's equation \cite{C}, 
Yor \cite{Y} studied the equation on the state space $ \bR $ given as 
\begin{align}
\eta_k = \xi_k + \{ \eta_{k-1} \} 
, \qquad k \in -\bN 
\label{eq: Yor}
\end{align}
for a general noise process $ \xi $, 
where $ \{ x \} $ stands for the fractional part of $ x $. 
He characterized the properties of the set of solutions 
in terms of the noise laws. 
Let us recall some of his results. 

Let $ \mu = (\mu_k:k \in -\bN) $ be a family of probability laws on $ \bR $. 
Define 
\begin{align}
\bZ_{\mu} = \cbra{ p \in \bZ : \pi_\mu(p) := 
\lim_{l \to -\infty } \prod_{k: k \le l} \absol{ \int_{\bR} \e^{2 \pi i p x} \mu_k(dx) } > 0 } . 
\label{eq: Z mu}
\end{align}
Note that $ \pi_{\mu}(p) = 1 $ if $ p \in \bZ_{\mu} $, 
while $ \pi_{\mu}(p)=0 $ otherwise. 
Then it follows (see \cite[Prop.3]{Y}) that 
$ \bZ_{\mu} $ is a subgroup of the additive group $ \bZ $. 
Now there exists a unique non-negative integer $ p_{\mu} $ such that 
$ \bZ_{\mu} = p_{\mu} \bZ $. 
The following theorem, 
which summarizes Prop.4, Thm.3, Thm.4 and Thm.5 of \cite{Y}, 
gives a complete answer to {\bf (Q1)}-{\bf (Q3)}: 

\begin{Thm}[{\cite{Y}}] \label{thm: Yor}
Let $ \{ (\eta_k),(\xi_k) \} $ denote any solution of \eqref{eq: Yor}. 
Then the following statements hold: 
\\ \quad {\rm (A)} 
Uniqueness in law holds if and only if $ p_{\mu}=0 $. 
In this case, 
it holds that $ \cF^{\eta}_{-\infty } $ is trivial and, 
for any $ k \in -\bN $, that the fractional part 
$ \{ \eta_k \} $ is uniform on $ [0,1) $, and that 
\begin{align}
\cF^{\eta}_k = \sigma(\{ \eta_k \}) \vee \cF^{\xi}_k 
\quad \text{a.s.} 
\label{}
\end{align}
where $ \sigma(\{ \eta_k \}) $ and $ \cF^{\xi}_k $ are independent; 
\\ \quad {\rm (B)} 
There exists a strong solution if and only if $ p_{\mu}=1 $. 
In this case, 
it holds, for any $ k \in -\bN $, that 
\begin{align}
\cF^{\eta}_k = \cF^{\eta}_{-\infty } \vee \cF^{\xi}_k 
\quad \text{a.s.} 
\label{}
\end{align}
where $ \cF^{\eta}_{-\infty } $ and $ \cF^{\xi}_k $ are independent; 
\\ \quad {\rm (C)} 
If $ p_{\mu} = 2,3,\ldots $, then 
it holds, for any $ k \in -\bN $, that 
\begin{align}
\{ p_{\mu} \eta_k \} \in \cF^{\eta}_{-\infty } \vee \cF^{\xi}_k \ \text{a.s.} 
, \quad 
[ p_{\mu} \eta_k ] \ \text{is uniform on } 
\{ 0,1,\ldots,p_{\mu}-1 \} , 
\label{}
\end{align}
where $ [x] $ stands for the integer part of $ x $, 
and that 
\begin{align}
\cF^{\eta}_k = \sigma([ p_{\mu} \eta_k ]) \vee \cF^{\eta}_{-\infty } \vee \cF^{\xi}_k 
\quad \text{a.s.} 
\label{}
\end{align}
where 
the three $ \sigma $-fields 
$ \sigma([ p_{\mu} \eta_k ]) $, 
$ \cF^{\eta}_{-\infty } $ and $ \cF^{\xi}_k $ are independent. 
\end{Thm}

Note that if $ \{ (\eta_k),(\xi_k) \} $ is a solution 
of equation \eqref{eq: Yor} taking values in $ \bR $, 
then the pair $ \{ (\{ \eta_k \}),(\{ \xi_k \}) \} $ is a solution 
of equation \eqref{eq: intro SE} taking values in $ \bR/\bZ $. 
Moreover, it has been proved in \cite[\S 9]{AUY} that 
certain properties, which are of interest to us, 
of the solutions of equation \eqref{eq: Yor} 
are equivalent to 
those of the solutions of equation \eqref{eq: intro SE}. 
Here $ \bZ_{\mu} = p_{\mu} \bZ $ 
corresponds to our $ H_{\mu} $ in relation \eqref{eq: H mu}. 

\begin{Rem}
{\rm (i)} 
Identity \eqref{eq: Z mu} shows 
how to compute the characteristic $ p_{\mu} $ 
from the noise law $ \mu $. 
We can characterize the subgroup $ H_{\mu} $ 
in terms of the noise law $ \mu $ 
completely in the case where $ G $ is commutative, 
but we do not know how to do this in the general case; 
see \cite[Thm.6.1]{AUY}. 
\\ \quad {\rm (ii)} 
Our Corollary \ref{thm: main5} gives more information than Yor's Theorem \ref{thm: Yor} 
in that the remote past $ \cF^{\eta}_{-\infty } $ is given explicitly 
as $ \sigma(V) $ in the cases {\rm (B)} and {\rm (C)}. 
\end{Rem}

\subsection{General lemmas}

Let us give several general lemmas 
concerning solutions of equation \eqref{eq: intro SE} taking values in compact groups. 

\begin{Lem}[{\cite[Lem.4.3]{AUY}}] \label{L2}
The following assertions hold: 
\\ \quad {\rm (i)} 
Let $ \{ (\eta^1_k),(\xi^1_k) \} $ 
and $ \{ (\eta^2_k),(\xi^2_k) \} $ be two solutions of \eqref{eq: intro SE}. 
Suppose that $ \eta^1_l \dist \eta^2_l $ for all $ l \in -\bN $. 
Then $ (\eta^1_k) \dist (\eta^2_k) $. 
\\ \quad {\rm (ii)} 
Let $ (\lambda_k) \subset \sP(G) $ 
which satisfies the convolution equation: 
\begin{align}
\lambda_k = \mu_k * \lambda_{k-1} 
, \quad k \in -\bN . 
\label{eq: conv eq2}
\end{align}
Then there exists a solution $ \{ (\eta_k),(\xi_k) \} $ 
such that each $ \eta_k $ has law $ \lambda_k $. 
\end{Lem}

\begin{proof}
(i) Let $ l \in -\bN $. 
For any $ k \ge l $, we have $ \eta^i_k = \xi^i_k \cdots \xi^i_{l+1} \eta^i_{l} $ 
for $ i=1,2 $. 
Hence we see that 
the joint laws of $ (\eta^i_{l},\eta^i_{l+1},\ldots,\eta^i_{0}) $ for $ i=1,2 $ 
coincide. 
This proves that $ (\eta^1_k) \dist (\eta^2_k) $. 

(ii) For any $ l \in -\bN $, we construct a family of random variables 
$ \{ \eta^{(l)}_k, \xi^{(l)}_k : k=l,\ldots,0 \} $ as follows: 
Let $ X,\xi_l,\ldots,\xi_0 $ be independent random variables such that 
$ X $ has law $ \lambda_{l-1} $ and each $ \xi_k $ has law $ \mu_k $. 
For $ k=l,\ldots,0 $, we define $ \eta^{(l)}_k = \xi_k \xi_{k-1} \cdots \xi_l X $. 
Then from the convolution equation \eqref{eq: conv eq2}, 
it follows easily that the family $ \{ \Pi^{(l)} : l \in -\bN \} $ 
of probability laws $ \Pi^{(l)} $ 
of $ \{ \eta^{(l)}_k, \xi^{(l)}_k : k=l,\ldots,0 \} $ 
is consistent. 
Thus, by Kolmogorov's extension theorem, 
we see that there exists a pair of processes $ \{ (\eta_k),(\xi_k) \} $ 
such that, for each $ l \in -\bN $, the law 
of $ \{ \eta_k, \xi_k : k=l,\ldots,0 \} $ is $ \Pi^{(l)} $. 
It is now easy to verify that the process $ \{ (\eta_k),(\xi_k) \} $ 
is as desired. 
\end{proof}

\begin{Thm}[\cite{Y},\cite{AUY}] \label{thm: unif sol}
There exists a unique uniform solution, 
i.e., a solution $ \{ (\eta^*_k),(\xi_k) \} $ 
such that each $ \eta^*_k $ is uniform on $ G $. 
Moreover, each $ \eta^*_k $ is independent of $ (\xi_k) $. 
\end{Thm}

\begin{proof}
Let $ \lambda_k = \omega_G $ for all $ k \in -\bN $. 
Then the sequence $ (\lambda_k) $ satisfies the convolution equation \eqref{eq: conv eq2}, 
and hence we obtain the desired conclusion by Lemma \ref{L2}. 
\end{proof}

\begin{Rem} \label{rem: unif sol}
A process $ (\eta_k) $ is called {\em stationary} 
if, for each $ n \in \bN $, 
the joint law of $ (\eta_{k},\eta_{k-1},\ldots,\eta_{k-n}) $ 
does not depend on $ k \in -\bN $. 
Since $ \xi_k = \eta^*_k (\eta^*_{k-1})^{-1} $, 
we see that, if the process $ (\eta^*_k) $ is stationary, 
then the noise $ (\xi_k) $ is identically distributed. 
In this case, the process $ (\eta^*_k) $ is stationary. 
See \cite{Takahashi} for the detailed discussion in this case. 
\end{Rem}

\begin{Lem} \label{thm: strong sol}
Let $ \{ (\eta_k),(\xi_k) \} $ be a solution. 
Then the following statements are equivalent: 
\\ \quad {\rm (i)} 
For all $ k \in -\bN $, $ \cF^{\eta}_k = \cF^{\xi}_k $ a.s. 
(i.e., this solution is strong); 
\\ \quad {\rm (ii)} 
There exists $ k \in -\bN $ such that $ \cF^{\eta}_k = \cF^{\xi}_k $ a.s.; 
\\ \quad {\rm (iii)} 
$ \cF^{\eta}_0 = \cF^{\xi}_0 $ a.s. 
\end{Lem}

\begin{proof}
It is obvious that {\rm (i)} implies {\rm (iii)} 
and that {\rm (iii)} implies {\rm (ii)}. 
Let us prove that {\rm (ii)} implies {\rm (i)}. 

Suppose that $ \cF^{\eta}_{k_0} = \cF^{\xi}_{k_0} $ a.s. for some $ k_0 \in -\bN $. 
For $ k \ge k_0+1 $, 
since $ \eta_k = \xi_k \cdots \xi_{k_0+1} \eta_{k_0} $ 
and since $ \eta_{k_0} \in \cF^{\xi}_{k_0} \subset \cF^{\xi}_k $ a.s., 
we have $ \eta_k \in \cF^{\xi}_k $. 
For $ k \le k_0-1 $, 
since $ \eta_k = (\xi_{k_0} \cdots \xi_{k+1})^{-1} \eta_{k_0} $, 
we have $ \eta_k \in \cF^{\xi}_{k_0} $. 
Since $ \cF^{\xi}_{k_0} = \cF^{\xi}_k \vee \cG $ 
with $ \cG = \sigma(\xi_{k+1},\ldots,\xi_{k_0}) $ 
which is assumed independent of $ \sigma(\eta_k) \vee \cF^{\xi}_k $, 
Lemma \ref{lem: independence} shows that 
$ 1_A(\eta_k) = E[1_A(\eta_k)|\cF^{\xi}_{k_0}] = E[1_A(\eta_k)|\cF^{\xi}_{k}] $ 
for all $ A \in \cB(G) $. 
This proves that $ \eta_k \in \cF^{\xi}_{k} $ a.s. 
Thus we obtain $ \cF^{\eta}_k \subset \cF^{\xi}_k $ a.s. for all $ k \in -\bN $. 
By identity \eqref{eq: noise subset eta}, 
we obtain $ \cF^{\eta}_k = \cF^{\xi}_k $ a.s. for all $ k \in -\bN $. 
The proof is now complete. 
\end{proof}

Recall that 
$ \sP_\mu $ is the set of the laws of $ (\eta_k) $ on $ G^{-\bN} $ 
for all possible solutions of equation \eqref{eq: intro SE}. 
Thus $ \sP_\mu $ is a subset of the compact convex set $ \sP(G^{-\bN}) $ 
where $ \sP(G^{-\bN}) $ is equipped with the topology of weak convergence. 
Moreover, by Markov property \eqref{eq: Markov}, we see that 
$ \sP_\mu $ is also compact and convex. 

\begin{Lem} \label{L3}
Let $ \{ (\eta_k),(\xi_k) \} $ be a solution. 
Then it is extremal, 
i.e., the law of $ (\eta_k) $ is an extremal point of $ \sP_\mu $, 
if and only if 
$ \cF^{\eta}_{-\infty } $ is trivial. 
\end{Lem}

The proof can be found in \cite[Thm.1]{Y} and \cite[Lem.1.2]{AUY}, 
but we give it for completeness of this paper. 

\begin{proof}
We may assume without loss of generality that $ (\eta_k) $ 
is the coordinate process on $ G^{-\bN} $. 

Suppose that $ \cF^{\eta}_{-\infty } $ is trivial. 
Suppose also that $ P $ can be represented as $ P=c P_1 + (1-c) P_2 $ 
for some $ P_1,P_2 \in \sP_\mu $ and $ 0<c<1 $. 
Then $ P_1 $ is absolutely continuous with respect to $ P $. 
By the Radon--Nikodym theorem, we see that 
there exists a non-negative functional $ D $ such that $ \d P_1 = D \d P $. 
Let $ Z $ be a non-negative functional. 
Since $ P $ and $ P_1 $ are solutions with the same noise, 
we have $ P[Z|\cF^{\eta}_k] = P_1[Z|\cF^{\eta}_k] $ 
$ P $-a.s. by Markov property \eqref{eq: Markov}. 
Hence we have 
\begin{align}
P[DZ] = P_1[Z] = P_1[P_1[Z|\cF^{\eta}_k]] = P[D P[Z|\cF^{\eta}_k] ]  
= P[P[D|\cF^{\eta}_k] Z] . 
\label{}
\end{align}
This shows that $ D \in \cF^{\eta}_{-\infty } $ $ P $-a.s. 
Since $ \cF^{\eta}_{-\infty } $ is $ P $-trivial, 
we have $ D=1 $ and $ P_1=P $. 
This proves that $ P $ is an extremal point of $ \sP_\mu $. 

Suppose that $ \cF^{\eta}_{-\infty } $ is not trivial. 
Then there exists a set $ A \in \cF^{\eta}_{-\infty } $ such that $ c:=P(A) \in (0,1) $. 
Hence $ P $ may be represented as $ P=c P_1 + (1-c) P_2 $ where 
$ \d P_1 = 1_A \d P /c $ and $ \d P_2 = 1_{A^c} \d P / (1-c) $. 
It is easy to see that $ P_1,P_2 \in \sP_\mu $, 
which shows that $ P $ is not an extremal point of $ \sP_\mu $. 
\end{proof}

\subsection{Results from Akahori--Uenishi--Yano \cite{AUY}}

Let us recall several results from Akahori et al.~\cite{AUY}. 

\begin{Thm}[{\cite{AUY}}] \label{L4}
Let $ \{ (\eta^1_k),(\xi^1_k) \} $ 
and $ \{ (\eta^2_k),(\xi^2_k) \} $ be two solutions of \eqref{eq: intro SE}. 
Suppose that they are extremal. 
Then there exists a deterministic element $ g \in G $ such that 
$ \{ (\eta^2_k),(\xi^2_k) \} \dist \{ (\eta^1_k g),(\xi^1_k) \} $. 
\end{Thm}

The proof of Theorem \ref{L4} can be found in \cite[Thm.1.3]{AUY}, 
which was based on a coupling method. So we omit the proof. 

\begin{Cor} \label{thm: Krein Milman}
Let $ \{ (\eta^0_k),(\xi_k) \} $ be an extremal solution. 
Then any solution $ \{ (\eta^1_k),(\xi^1_k) \} $ may be represented as 
\begin{align}
\{ (\eta^1_k),(\xi^1_k) \} \dist \{ (\eta^0_k V),(\xi_k) \} 
\label{}
\end{align}
for some $ G $-valued random variable $ V $ 
independent of $ \{ (\eta^0_k),(\xi_k) \} $. 
\end{Cor}

\begin{proof}
From Theorem \ref{L4}, it follows that 
the laws $ P^{(\eta^0_k g)} $ of $ (\eta^0_k g) $ for $ g \in G $ 
exhaust all extremal points of $ \sP_\mu $. 
By the Krein--Milman theorem (see, e.g., \cite{Ph}), we see that 
the law $ P^{(\eta^1_k)} $ of $ (\eta^1_k) $ may be represented as 
\begin{align}
P^{(\eta^1_k)}(\cdot) = \int_G P^{(\eta^0_k g)}(\cdot) \nu(\d g) 
\label{}
\end{align}
for some probability law $ \nu $ on $ G $. 
Then we have $ (\eta^1_k) \dist (\eta^0_k V) $ 
for some $ G $-valued random variable $ V $ 
independent of $ (\eta^0_k) $. 
By equation \eqref{eq: noise subset eta}, 
we complete the proof. 
\end{proof}

Akahori et al.~\cite{AUY} partially generalized Yor's Theorem \ref{thm: Yor}. 
To summarize in the framework of groups, 
we may say that Yor's study \cite{Y} was based on the {\em Pontryagin duality} 
between the (locally) compact group $ \bR / \bZ $ and 
the class of all characters on $ \bR / \bZ $, 
while Akahori et al.~\cite{AUY} was based on the {\em Tannaka duality} 
between a compact group $ G $ and the class of all unitary representations on $ G $. 

Let $ \{ (\eta^0_k),(\xi_k) \} $ be an extremal solution. 
Define 
\begin{align}
H^{\rm iso}_\mu = \{ g \in G : (\eta^0_k g) \dist (\eta^0_k) \}. 
\label{}
\end{align}
Let $\fG$ denote the set of all unitary representations $ \rho $ of $ G $ 
on a finite dimensional linear space. 
Define 
\begin{align}
\fH^{\rm strong}_\mu 
= \cbra{ \rho \in \fG : \rho(\eta^0_k) \in \cF^\xi_k \ \text{a.s. for all $k \in -\bN$} } 
\label{}
\end{align}
and 
\begin{align}
H^{\rm strong}_\mu = \cbra{ g \in G : \rho(g) = {\rm id} 
\ \text{for every $\rho \in \fH^{\rm strong}_\mu$} } . 
\label{}
\end{align}

\begin{Thm}[{\cite[Thm.1.6]{AUY}}] \label{T3}
The following statements hold: 
\\ \quad {\rm (A)} 
Uniqueness in law holds if and only if 
$ H^{\rm iso}_\mu = G $; 
\\ \quad {\rm (B)} 
There exists a strong solution if and only if 
$ H^{\rm strong}_\mu = \{ {\rm unit} \} $. 
\end{Thm}

By virtue of our Theorem \ref{thm: main1}, 
we obtain the following theorem: 

\begin{Thm} \label{thm: main6}
Let $ H_{\mu} $ 
and $ \{ (\eta^0_k),(\xi_k) \} $ 
as in Theorem \ref{thm: main1}. 
Let $ H^{\rm iso}_\mu $ and $ H^{\rm strong}_\mu $ 
be associated with $ \{ (\eta^0_k),(\xi_k) \} $. 
Then it holds that 
\begin{align}
H^{\rm iso}_\mu = H_\mu 
, \quad 
H^{\rm strong}_\mu = \bigcup_{g \in G} g H_{\mu} g^{-1} . 
\label{}
\end{align}
In other words, $ H^{\rm strong}_\mu $ 
is the smallest normal subgroup containing $ H_{\mu} $. 
\end{Thm}

The proof of Theorem \ref{thm: main6} will be given in Subsection \ref{sec: AUY proof}.

\section{Proofs of main theorems} \label{sec: infin prod}

We prove our main theorems in the following order: 
Theorem \ref{thm: main0}, Theorem \ref{thm: main1.2}, Theorem \ref{thm: main1.4}, 
Theorem \ref{thm: main1}, and Theorem \ref{thm: main3}.

\subsection{General principle of Kloss--Tortrat--Csisz{\'a}r}

Limit laws of infinite products of random variables on compact groups 
have been first studied by Kawada--It\^o \cite{KI}. 
After that, Kloss \cite{K} discovered a general principle of infinite products, 
which was generalized to locally compact groups 
by Tortrat \cite{T} and by Csisz{\'a}r \cite{Csi} independently. 
Let us recall some results from Csisz{\'a}r \cite{Csi}. 
For some discussions in the case of locally compact semigroups, 
see, e.g., Mukherjea--Tserpes \cite{MT}. 
For basic notations and facts about probability laws on compact groups, 
see, e.g., standard textbooks \cite{P}, \cite{H} and \cite{G}. 

The following theorems are taken from Csisz\'ar \cite{Csi}, 
where he called them {\em Kloss's general principle}: 

\begin{Thm}[{\cite[Thm.3.1]{Csi}}] \label{C2}
Let $ (\xi_n: n \in \bN ) $ be a sequence of independent $ G $-valued random variables.
Then there exists a sequence $ (\alpha_m: m \in \bN) $ 
of deterministic elements of $ G $ such that, 
for any $ n \in \bN $, the product 
$ \xi_n \xi_{n+1} \cdots \xi_{m} \alpha _{m} $ 
converge in law as $ m \to \infty $. 
\end{Thm}

\begin{Thm}[{\cite[Thm.3.2]{Csi}}] \label{C3}
Let $ (\xi_n: n \in \bN ) $ be a sequence of independent $ G $-valued random variables.
Assume that, for each $ n \in \bN $, 
the product $ \xi_n \xi_{n+1} \cdots \xi_{m} $ 
converges in law as $ m \to \infty $ 
to some $ G $-valued random variable $ \eta_n $. 
Then there exists a unique compact subgroup $ H $ of $ G $ such that 
the following statements hold: 
\\ \quad {\rm (i)} 
For each $ n \in \bN $ and $ h \in H $, $ \eta_n h \dist \eta_n $; 
\\ \quad {\rm (ii)} 
For each $ n \in \bN $, $ \xi_n \xi_{n+1} \cdots \xi_{m} H $ 
converges a.s. in $ G/H $ as $ m \to \infty $. 
\\
In this case, it holds that 
\\ \quad {\rm (iii)} 
$ \eta_n \cdist \omega_H $ as $ n \to \infty $. 
\end{Thm}

\subsection{Infinite products of independent random variables} \label{sec: ex sol}

For any probability laws $ (\mu_k) $ on $ G $ 
and any $ k,l \in -\bN $ with $ k \ge l $, we write 
$ \mu_{k,l} := \mu_k * \mu_{k-1} * \cdots * \mu_{l} $; 
for instance, $ \mu_{k,k} = \mu_k $, 
$ \mu_{k,k-1} = \mu_k * \mu_{k-1} $, and so on. 
For any $ G $-valued random variables $ (\xi_k) $ 
and any $ k,l \in -\bN $ with $ k \ge l $, we write 
$ \xi_{k,l} = \xi_k \xi_{k-1} \cdots \xi_l $; 
for instance, $ \xi_{k,k} = \xi_k $, $ \xi_{k,k-1} = \xi_k \xi_{k-1} $, and so on. 

Let us prove Theorem \ref{thm: main0}. 

\noindent
{\it Proof of Theorem \ref{thm: main0}.} 
By Theorem \ref{C2}, 
there exist a sequence $ (\alpha _l) $ of deterministic elements of $ G $ 
and a sequence $ (\lambda_k) $ of probability laws on $ G $ 
such that, 
for any $ k \in -\bN $, 
$ \mu_{k,l} * \delta_{\alpha _{l}} \to \lambda_k $ 
as $ l \to -\infty $. 
This shows {\bf (I1)}. 

It is obvious that $ (\lambda_k) $ 
satisfies the convolution equation \eqref{eq: conv eq2}. 
Hence, by Lemma \ref{L2}, there exists a solution $ \{ (\eta^0_k),(\xi_k) \} $ 
such that, for any $ k \in -\bN $, it holds that 
$ \eta^0_k $ has law $ \lambda_k $ and that 
\begin{align}
\xi_{k,l} \alpha _{l} 
\cdist \eta^0_k 
\quad \text{as $ l \to -\infty $}. 
\label{eq: pre E1}
\end{align}
Set $ \tilde{\xi}_k = \alpha _{k+1}^{-1} \xi_k \alpha _k $ for $ k \in -\bN $. 
Then we see, for any $ k \in -\bN $, that 
\begin{align}
\tilde{\xi}_{k,l} 
\cdist \alpha _{k+1}^{-1} \eta^0_k 
\quad \text{as $ l \to -\infty $} . 
\label{}
\end{align}
Then Theorem \ref{C3} shows that 
there exists a compact subgroup $ H $ of $ G $ such that 
\\ \quad {\rm (i)} 
$ \eta^0_k h \dist \eta^0_k $ 
for each $ h \in H $ and each $ k \in -\bN $; 
\\ \quad {\rm (ii)} 
$ \tilde{\xi}_{k,l} H $ 
converges a.s. in $ G/H $ as $ l \to -\infty $ 
for each $ k \in -\bN $; 
\\ \quad {\rm (iii)} 
$ \alpha _l^{-1} \eta^0_{l-1} \cdist \omega_H $ as $ l \to -\infty $. 
\\
The statement {\rm (ii)} implies 
\\ \quad {\rm (ii$ ' $)} 
$ \xi_{k,l} \alpha _{l} H $ 
converges a.s. in $ G/H $ as $ l \to -\infty $ 
for each $ k \in -\bN $. 

The statements {\rm (iii)}, {\rm (i)} and {\rm (ii$ ' $)} 
prove {\bf (I2)}, {\bf (I3)} and {\bf (I4)}, respectively. 

Suppose that $ (\tilde{\lambda}_k) $, $ (\tilde{\alpha }_k) $ and $ \tilde{H} $ 
also satisfy {\bf (I1)} and {\bf (I2)} 
and let $ g $ be an accumulation point such that 
$ \alpha _{l}^{-1} \tilde{\alpha }_{l} \to g $ 
along a subsequence $ l=l(j) $. 

By {\bf (I1)}, we have 
$ \mu_{k,l} * \delta_{\alpha _{l}} \to \lambda_k $ 
and $ \mu_{k,l} * \delta_{\tilde{\alpha} _{l}} \to \tilde{\lambda}_k $ 
as $ l \to -\infty $ 
for each $ k \in -\bN $. 
Taking the limit 
in both sides of the identity 
\begin{align}
\mu_{k,l} * \delta_{\tilde{\alpha} _{l}} 
= \mu_{k,l} * \delta_{\alpha _{l}} * \delta_{ \alpha _{l}^{-1} \tilde{\alpha} _{l} } 
\label{}
\end{align}
along the subsequence $ l=l(j) $, 
we have $ \tilde{\lambda}_k = \lambda_k * \delta_g $ for all $ k \in -\bN $. 

By {\bf (I2)}, we have 
$ \delta_{\alpha _{l}^{-1}} * \lambda_{l-1} \to \omega_H $ 
and $ \delta_{\tilde{\alpha} _{l}^{-1}} * \tilde{\lambda}_{l-1} \to \omega_{\tilde{H}} $ 
as $ l \to -\infty $. 
Taking the limit 
in both sides of the identity 
\begin{align}
\delta_{\tilde{\alpha} _{l}^{-1}} * \tilde{\lambda}_{l-1} 
= \delta_{ \rbra{ \alpha _{l}^{-1} \tilde{\alpha} _{l} }^{-1}} 
* \delta_{\alpha _{l}^{-1}} * \lambda_{l-1} * \delta_g 
\label{}
\end{align}
along the subsequence $ l=l(j) $, 
we obtain $ \omega_{\tilde{H}} = \delta_{g^{-1}} * \omega_H * \delta_g $, 
which proves that $ \tilde{H} = g^{-1} H g $. 
Therefore the proof of Theorem \ref{thm: main0} is complete. 
\qed

Let us prove Theorem \ref{thm: main1.2}. 

\noindent
{\it Proof of Theorem \ref{thm: main1.2}.} 
{\rm (i)} 
Let $ \{ (\eta^0_k),(\xi_k) \} $ be a solution 
satisfying {\bf (E1)} with some $ (\alpha _l) \subset G $. 
Let us prove that the solution $ \{ (\eta^0_k),(\xi_k) \} $ is extremal. 

Let $ \{ (\eta_k),(\xi_k) \} $ be an arbitrary solution. 
Here we denote the noise process by the same notation 
without any confusion. 
Note that
\begin{align}
\eta_k 
= \xi_{k,l} \eta_{l-1} 
= (\xi_{k,l} \alpha_{l}) (\alpha_{l}^{-1} \eta_{l-1}) . 
\label{**}
\end{align}
Since $ \sP(G) $ is compact, there exists a subsequence $ l=l(j) $ such that
\begin{align}
\alpha_{l}^{-1} \eta_{l-1} \cdist V 
\quad \text{along $ l=l(j) $}
\label{eq: pre E1 *}
\end{align}
for some $ G $-valued random variable $ V $, 
which we may take to be independent of $ \{ (\eta^0_k),(\xi_k) \} $. 

Now we take the limit in \eqref{**} along the subsequence $ l=l(j) $. 
Note that, since $ \{ (\eta_k),(\xi_k) \} $ is a solution, 
we see that $ \alpha^{-1}_{l} \eta_{l-1} $ is independent 
of $ \xi_{k,l} \alpha_{l} $. 
By \eqref{eq: pre E1} and \eqref{eq: pre E1 *}, 
we see that 
\begin{align}
\rbra{ \xi_{k,l} \alpha_{l} \, , \, \alpha^{-1}_{l} \eta_{l-1} } 
\cdist (\eta^0_k,V) 
\quad \text{along $ l=l(j) $} . 
\end{align}
Taking the limit in \eqref{**} along this subsequence, 
we have $ \eta_k \dist \eta^0_k V $ for each $ k \in -\bN $. 
By Lemma \ref{L2}, we have $ (\eta_k) \dist (\eta^0_k V) $. 
This proves that the solution $ \{ (\eta^0_k),(\xi_k) \} $ is extremal. 

(ii) 
Let $ \{ (\eta^0_k),(\xi_k) \} $ be a solution 
which satisfies {\bf (E1)}, {\bf (E3)} and {\bf (E4)} 
for some sequence $ (\alpha _l) \subset G $ 
and some compact subgroup $ H $ of $ G $. 
Applying Theorem \ref{C3} 
to $ \tilde{\xi}_k = \alpha _{k+1}^{-1} \xi_k \alpha _k $, 
we see that {\bf (E2)} holds with these $ (\alpha _l) $ and $ H $. 

The proof of Theorem \ref{thm: main1.2} is therefore complete. 
\qed

\subsection{Construction of an extremal solution} \label{sec: construction}

Let us prove Theorem \ref{thm: main1.4}. 

\noindent
{\it Proof of Theorem \ref{thm: main1.4}.} 
Let $ \xi_k $'s be independent random variables such that 
each $ \xi_k $ has law $ \mu_k $. 
Let $ U_0 $ be a $ G $-valued random variable 
which is independent of $ (\xi_k) $ and is uniform on $ H $. 

By Theorem \ref{thm: main0}, we see that, for any $ k \in -\bN $, the limit 
\begin{align}
\Phi_k := \lim_{l \to -\infty } \xi_{k,l} \alpha _{l} H 
\label{}
\end{align}
converges in $ G/H $ a.s. as $ l \to -\infty $. 
Hence, for any fixed $ k \in -\bN $, 
we may define $ \phi_k $ as $ \phi_k = \vs(\Phi_k) $. 
It is by definition 
that $ \phi_k $ is a.s. measurable with respect to $ \cF^{\xi}_k $ and that 
\begin{align}
\xi_{k,l} \alpha _{l} H \asto \phi_k H 
\quad \text{as $ l \to -\infty $}. 
\label{**2}
\end{align}
Now it is obvious that 
\begin{align}
\phi_k H = \xi_k \phi_{k-1} H 
\quad \text{a.s. for all $ k \in -\bN $}. 
\label{*2}
\end{align}

For $ k \in -\bN $, we define 
\begin{align}
U_k = \phi_k^{-1} \xi_{0,k+1}^{-1} \phi_0 U_0 . 
\label{*2**}
\end{align}
Note that 
\begin{align}
h_k := \phi_k^{-1} \xi_{0,k+1}^{-1} \phi_0 
= \rbra{ \phi_k^{-1} \xi_{k+1}^{-1} \phi_{k+1} } 
\rbra{ \phi_{k+1}^{-1} \xi_{k+2}^{-1} \phi_{k+2} } 
\cdots 
\rbra{ \phi_{-1}^{-1} \xi_{0}^{-1} \phi_{0} } , 
\label{}
\end{align}
which belongs to $ H $ a.s. by \eqref{*2}. Hence we see that 
\begin{align}
\rbra{ (\xi_k),U_k } 
= \rbra{ (\xi_k), h_k U_0 } 
\dist \rbra{ (\xi_k),U_0 } ; 
\label{*2***}
\end{align}
in fact, since $ h_k \in \cF^{\xi}_0 $, we have, 
for any bounded measurable function $ f $ on $ G $, 
\begin{align}
E[f(h_k U_0)|\cF^{\xi}_0] = 
\int_H f(h_k h) \omega_H(\d h) 
= \int_H f(h) \omega_H(\d h) 
= E[f(U_0)] . 
\label{}
\end{align}
Now we see by \eqref{*2***} that 
$ U_k $ is independent of $ \cF^{\xi}_0 $ and is uniform on $ H $.

We define 
\begin{align}
\eta^0_k = \phi_k U_k 
, \quad k \in -\bN . 
\label{*2*}
\end{align}
By \eqref{*2**}, we have $ \eta^0_k = \xi_k \eta^0_{k-1} $ a.s. 
for each $ k \in -\bN $. 
Let us prove that 
each $ \xi_k $ is independent of $ \cF^{\eta^0}_{k-1} $. 
Let $ k>l $ and let $ f_k,f_{k-1},\ldots,f_l $ be non-negative measurable functions on $ G $. 
Then we have 
\begin{align}
& E[ f_k(\xi_k) f_{k-1}(\eta^0_{k-1}) f_{k-2}(\eta^0_{k-2}) \cdots 
f_l(\eta^0_l) ] 
\label{} \\
=& E[ f_k(\xi_k) f_{k-1}(\eta^0_{k-1}) 
f_{k-2}((\xi_{k-1,k-1})^{-1} \eta^0_{k-1}) 
\cdots 
f_l((\xi_{k-1,l+1})^{-1} \eta^0_{k-1}) ] 
\label{} \\
=& E[ f_k(\xi_k) f_{k-1}(\phi_{k-1} U_{k-1}) 
f_{k-2}(\psi_{k-1} U_{k-1}) 
\cdots 
f_l(\psi_{l+1} U_{k-1}) ] 
\label{** eq1}
\end{align}
where $ \psi_j = (\xi_{k-1,j})^{-1} \phi_{k-1} $ for $ j=k-1,\ldots,l+1 $. 
By \eqref{*2***}, we obtain 
\begin{align}
\text{\eqref{** eq1}} 
=& E[ f_k(\xi_k) f_{k-1}(\phi_{k-1} U_0) 
f_{k-2}(\psi_{k-1} U_0) 
\cdots 
f_l(\psi_{l+1} U_0) ] . 
\label{** eq2}
\end{align}
Since $ \xi_k $ is independent of $ \cF^{\xi}_{k-1} \vee \sigma(U_0) $, 
we obtain 
\begin{align}
\text{\eqref{** eq2}} 
=& E[ f_k(\xi_k) ] 
E[ f_{k-1}(\phi_{k-1} U_0) 
f_{k-2}(\psi_{k-1} U_0) 
\cdots 
f_l(\psi_{l+1} U_0) ] . 
\label{***}
\end{align}
This proves that $ \xi_k $ is independent of 
$ \sigma(\eta^0_{k-1},\ldots,\eta^0_l) $, 
and hence of $ \cF^{\eta^0}_{k-1} $ by a monotone class argument. 
Therefore, we see that 
$ \{ (\eta^0_k),(\xi_k) \} $ is a solution. 

Let $ k \in -\bN $. 
By \eqref{**2}, we have 
\begin{align}
\xi_{k,l} \alpha _{l} H \asto \eta^0_k H 
\quad \text{as $ l \to -\infty $}. 
\label{**2+}
\end{align}
By \eqref{**2+} and by definition $ \eta^0_k = \phi_k U_k $, we have 
\begin{align}
\xi_{k,l} \alpha _{l} U_k \cdist \eta^0_k U_k \dist \eta^0_k 
\quad \text{as $ l \to -\infty $}. 
\label{}
\end{align}
On the other hand, by {\bf (I1)}, we see that 
the law of $ \xi_{k,l} \alpha_l U_k $ converges to 
$ \lambda_k * \omega_H $, 
which is equal to $ \lambda_k $ by {\bf (I3)}. 
Thus we conclude that $ \eta^0_k $ has law $ \lambda_k $. 

By (i) of Theorem \ref{thm: main1.2}, 
we see that the solution $ \{ (\eta^0_k),(\xi_k) \} $ is extremal. 
The proof of Theorem \ref{thm: main1.4} is therefore complete. 
\qed

\subsection{Characterization of extremal solutions} \label{sec: char ex solution}

Now we prove Theorem \ref{thm: main1}. 

\noindent
{\it Proof of Theorem \ref{thm: main1}.} 
Let $ \{ (\eta^0_k),(\xi_k) \} $, $ (\alpha _l) $ and $ H $ 
be as are given in the proof of Theorem \ref{thm: main0}. 
By (i) of Theorem \ref{thm: main1.2}, 
we see that $ \{ (\eta^0_k),(\xi_k) \} $ is an extremal solution 
satisfying {\bf (E1)}, {\bf (E2)} and {\bf (E3)} 
with $ (\alpha _l) $ and $ H $. 
We see by Theorem \ref{thm: main1.4} that 
$ \{ (\eta^0_k),(\xi_k) \} $ also satisfies {\bf (E4)}. 
Hence we see that 
this particular extremal solution $ \{ (\eta^0_k),(\xi_k) \} $ 
satisfies {\bf (E1)}-{\bf (E4)}. 
Since the general case follows immediately by Theorem \ref{L4}, 
we have now proved the former half of Theorem \ref{thm: main1}. 
The latter half of Theorem \ref{thm: main1} 
is immediate from that of Theorem \ref{thm: main0}. 
The proof of \ref{thm: main1} is therefore complete. 
\qed

\subsection{Characterization of $ H_{\mu} $} \label{sec: char subgr}

Let us prove Corollary \ref{thm: main1.5}. 
Before doing this, we prove the following lemma. 

\begin{Lem} \label{lem: conv cosets}
Let $ H $ and $ K $ be compact subgroups of $ G $. 
Let $ g,g_1,g_2,\ldots $ be elements of $ G $. 
Then the following statements hold: 
\\ \quad {\rm (i)} 
If $ H \subset K $ and if $ g_n H \to g H $, then $ g_n K \to g K $; 
\\ \quad {\rm (ii)} 
If $ g_n H \to g H $ and if $ g_n K \to g K $, then 
$ g_n (H \cap K) \to g (H \cap K) $. 
\end{Lem}

\noindent
{\it Proof of Lemma \ref{lem: conv cosets}.} 
(i) 
Let $ \pi_H $ and $ \pi_K $ denote the natural projections of $ G $ 
onto $ G/H $ and $ G/K $, respectively. 
Since $ H \subset K $, there exists a mapping $ \pi_{H,K}:G/H \to G/K $ such that 
$ \pi_K = \pi_{H,K} \circ \pi_H $. 
Then it is immediate that $ \pi_{H,K} $ is continuous. 
Hence we see that 
$ g K = \pi_{H,K}(gH) = \lim_n \pi_{H,K}(g_nH) = \lim_n g_n K $. 

(ii) 
Let $ \tilde{g}(H \cap K) $ be an accumulation point of $ \{ g_n (H \cap K) \} $, 
which exists by compactness of $ G/(H \cap K) $. 
We may take a subsequence $ n=n(j) $ 
such that $ g_{n(j)} (H \cap K) \to \tilde{g}(H \cap K) $. 
By (i), we see that 
$ g_{n(j)} H \to \tilde{g} H $ and $ g_{n(j)} K \to \tilde{g} K $, 
which implies that $ \tilde{g} H = gH $ and $ \tilde{g} K = gK $. 
This shows that $ \tilde{g} (H \cap K) = g (H \cap K) $. 
Thus we obtain $ g_n (H \cap K) \to g (H \cap K) $. 
\qed

\noindent
{\it Proof of Corollary \ref{thm: main1.5}.} 
{\rm (i)} 
Set $ H = \{ h \in G : (\eta^0_k h) \dist (\eta^0_k) \} $. 
Then it is obvious that $ H $ is a compact subgroup of $ G $ 
and contains $ H_{\mu} $. 
Then it is obvious that {\bf (E3)} holds with $ H $. 
By (i) of Lemma \ref{lem: conv cosets}, 
we see that {\bf (E4)} holds with $ (\alpha _l) $ and $ H $. 
Then, by (ii) of Theorem \ref{thm: main1.2}, 
it also satisfies {\bf (E2)} with $ (\alpha _l) $ and $ H $. 
This proves that $ H = H_{\mu} $. 

{\rm (ii)} 
Suppose that $ H $ is a compact subgroup such that 
\begin{align}
\xi_{k,l} \alpha _{l} H \asto \eta^0_k H 
\quad \text{as $ l \to -\infty $ for all $ k \in -\bN $}. 
\label{}
\end{align}
Set $ \tilde{H} = H \cap H_{\mu} $. 
By (ii) of Lemma \ref{lem: conv cosets}, 
we see that 
\begin{align}
\xi_{k,l} \alpha _{l} \tilde{H} \asto \eta^0_k \tilde{H} 
\quad \text{as $ l \to -\infty $ for all $ k \in -\bN $}. 
\label{}
\end{align}
Hence {\bf (E3)} and {\bf (E4)} hold with $ (\alpha _l) $ and $ \tilde{H} $. 
In the same way as above, we obtain $ \tilde{H} = H_{\mu} $, 
which implies that $ H \supset H_{\mu} $. 
\qed

\noindent
{\it Proof of Corollary \ref{thm: main2}.} 
This is obvious 
from Theorem \ref{thm: main1}, Corollary \ref{thm: main1.5} and 
Theorem \ref{L4}. 
\qed

\subsection{Complementation formulae} \label{sec: compl}

In this section, 
we let $ (\lambda_k) $, $ (\alpha _l) $ and $ H $ be as in Theorem \ref{thm: main0} 
and let $ \vs(\cdot) $ and $ \vh(\cdot) $ as in Subsection \ref{sec: factor and compl}. 
For $ g \in G $, we write 
$ \vs(g) $ simply for $ \vs(gH) $. 

Now we prove Theorem \ref{thm: main3}. 

\noindent
{\it Proof of Theorem \ref{thm: main3}.} 
Let $ \{ (\eta_k),(\xi_k) \} $ be any solution. 
Let $ U_0 $ be a $ G $-valued random variable 
which is independent of $ \{ (\eta_k),(\xi_k) \} $ 
and define $ (U_k) $ and $ (\eta^0_k) $ as given in Theorem \ref{thm: main1.4}. 
Since $ \{ (\eta^0_k),(\xi_k) \} $ is an extremal solution, 
there exists a $ G $-valued random variable $ V $ such that 
$ \{ (\eta_k),(\xi_k) \} \dist \{ (\eta^0_k V),(\xi_k) \} $. 
Noting that 
\begin{align}
\eta^0_k V = \xi_{0,k+1}^{-1} \phi_0 U_0 \vh(V^{-1})^{-1} \vs(V^{-1})^{-1} , 
\label{}
\end{align}
and that $ U_0 $ is independent of $ \cF^{\xi}_0 \vee \sigma(V) $, we have 
\begin{align}
\{ (\eta^0_k V),(\xi_k) \} \dist \{ (\eta^0_k \vs(V^{-1})^{-1}),(\xi_k) \} . 
\label{}
\end{align}
Thus we may assume without loss of generality that 
$ V = \vs(V^{-1})^{-1} $. 
For simplicity, let us write 
\begin{align}
\eta_k = \eta^0_k V = \phi_k U_k V . 
\label{eq: pf of main3}
\end{align}
Now it is obvious that Claims (i) and (ii) hold 
and that the three $ \sigma $-fields $ \sigma(U_k) $, $ \sigma(V) $ and $ \cF^{\xi}_k $ 
are independent. 

Let $ k \in -\bN $ be fixed. 
Since $ \eta_k = \phi_k U_k V $, we have $ \eta_k^{-1} \phi_k = V^{-1} U_k $. 
This shows that 
\begin{align}
V^{-1} = \vs(V^{-1}) = \vs(\eta_k^{-1} \phi_k) . 
\label{}
\end{align}
Since $ k $ is arbitrary, we obtain \eqref{eq: main4 02} 
and $ \sigma(V) \subset \cF^{\eta}_{-\infty } $. 
By \eqref{eq: pf of main3}, 
we obtain \eqref{eq: main3 sigma fields}. 

By Lemma \ref{lem: interchange}, we obtain 
\begin{align}
\cF^{\eta}_{-\infty } = \bigcap_{k \in -\bN} \cF^{\eta}_k 
\subset \bigcap_{k \in -\bN} \rbra{ \cF^{\eta^0}_k \vee \sigma(V) } 
= \rbra{ \bigcap_{k \in -\bN} \cF^{\eta^0}_k } \vee \sigma(V) 
= \sigma(V) 
\label{}
\end{align}
where we have used the fact that $ \cF^{\eta^0}_{-\infty } $ is trivial. 
Thus we obtain Claim (iii). 

Therefore the proof is complete. 
\qed

\subsection{Characteristic subgroups $ H^{\rm iso}_\mu $ 
and $ H^{\rm strong}_\mu $} \label{sec: AUY proof}

Now we prove Theorem \ref{thm: main6}. 

\noindent
{\it Proof of Theorem \ref{thm: main6}.} 
{\rm (i)} 
This is obvious by (i) of Corollary \ref{thm: main1.5}.

{\rm (ii)} 
Let us simply write $ H $ for $ H_{\mu} $. 
Set $ N_H = \bigcup_{g \in G} g H g^{-1} $. 
Since $ G $ and $ H $ are compact, 
we see that $ N_H $ is also compact. 
In fact, if $ g_n h_n g_n^{-1} \to f \in G $, 
then there exists a subsequence $ n(m) $ such that 
$ g_{n(m)} \to g \in G $ and $ h_{n(m)} \to h \in H $, 
and hence $ f = ghg^{-1} \in N_H $. 

Let us prove that $ H^{\rm strong}_{\mu} = N_H $. 

Let $ k \in -\bN $ be fixed. 
By the proof of Theorem \ref{thm: main1}, 
we may represent $ \eta^0_k $ as $ \eta^0_k = \phi_k U_k $ 
where 
$ \phi_k $ is measurable with respect to $ \cF^{\xi}_k $ 
and $ U_k $ is independent of $ (\xi_k) $ and is uniform on $ H $. 
Then, for any $ \rho \in \fH^{\rm strong}_\mu $, we have 
\begin{align}
\rho(U_k) 
= \rho(\phi_k)^{-1} \rho(\phi_k U_k) = \rho(\phi_k)^{-1} \rho(\eta^0_k) 
\in \cF^{\xi}_k 
\quad \text{a.s.}. 
\label{}
\end{align}
But, since $ \rho(U_k) $ is independent of $ (\xi_k) $, we see that 
$ \rho(U_k) $ is constant a.s. 
That is, $ \rho(h) $ is constant for $ \omega_H $-a.e.~$ h $. 
By continuity of $ \rho $, we have 
$ \rho $=id. on $ H $, which implies that 
$ \rho $=id. on $ N_H $. 
Now we obtain 
\begin{align}
\fH^{\rm strong}_\mu 
=& \cbra{ \rho : \rho(h) = \text{id. for every $h \in N_{H} $} } 
\label{} \\
=& \cbra{ \rho=(\rho_1,\rho_2) : 
\rho_1 = \mbox{id. of } N_{H}, \mbox{ $\rho_2$ 
is unitary repre. on $ G/N_H $} }. 
\label{}
\end{align}
Since $ N_H $ is a compact normal subgroup, 
the quotient $ G/N_H $ is again a compact group. 
Hence the stabilizer $ H^{\rm strong}_\mu $ of $ \fH^{\rm strong}_\mu $ 
is nothing else but $ N_H $. 

The proof of Theorem \ref{thm: main6} is therefore complete.
\qed


{\bf Acknowledgements.} 
The authors would like to thank Professors Marc Yor and Jean-Paul Thouvenot 
for their fruitful comments.


\def\cprime{$'$} \def\cprime{$'$}


\end{document}